\documentclass[a4paper,12pt]{article}
\usepackage[margin=28mm]{geometry}
\usepackage{amsmath,amssymb,mathtools}
\usepackage[T1]{fontenc}
\usepackage{lmodern}
\usepackage{microtype}
\usepackage{xcolor}
\usepackage{cite}
\usepackage[hidelinks,hypertexnames=false]{hyperref}
\usepackage[nameinlink,noabbrev]{cleveref}
\usepackage{xspace}
\DeclareRobustCommand{\eg}{\textup{eg}\xspace}
\DeclareRobustCommand{\ie}{\textup{ie}\xspace}

\newcommand{\R}{\mathbb R}
\newcommand{\ind}{\delta}
\newcommand{\prox}{\operatorname{prox}}
\newcommand{\diag}{\operatorname{diag}}
\newcommand{\ip}[2]{\langle #1,#2\rangle}
\newcommand{\res}{\mathcal R}
\newtheorem{proposition}{Proposition}
\AddToHook{env/proposition/begin}{\samepage}

\begin{document}
\title{A Unified View of Operator Splitting \\ and Interior-Point Methods}
\author{Gabriel Berk Pereira and Paul J. Goulart\\[0.4em]\small Department of Engineering Science, University of Oxford\\\small Oxford, United Kingdom}
\date{}
\maketitle

\begin{abstract}
First-order splitting and interior-point methods offer contrasting computational trade-offs for convex conic quadratic programs, yet their linear systems look remarkably similar. We make this connection precise through the Moreau slack, an unconstrained variable from which multipliers and (ordinary) slacks satisfying exact or barrier-perturbed complementarity are recovered. The remaining feasibility equations define a unified residual whose roots are KKT points at zero barrier parameter and central-path points at positive barrier parameter. Applying splitting to this residual yields projection methods at zero barrier parameter and barrier splitting at positive parameter; Newton's method in the latter setting gives an implicit interior-point method. At fixed positive barrier parameter, matching the splitting metric to current barrier curvature equalises the primal variable corrections under the same regularisation and leaves a factor of two between the Moreau-slack corrections. With metrics frozen at a barrier solution, local analysis identifies when doubling improves (local) convergence and when the primal solve remains the bottleneck. Metric choices that simplify linear solves (\eg make them diagonal) become Newton regularisations, exposing a shared convergence cost associated with computational convenience. Finally, a computable defect in the regularised Newton equations gives an accuracy-based criterion for reusing a splitting factorisation.
\end{abstract}

\begingroup\small\renewcommand{\contentsname}{\normalsize Contents}\tableofcontents\endgroup

\section{Introduction}

Throughout this paper we consider the following convex conic program:
\begin{equation}
	\begin{aligned}
	\min_x\quad &f(x):=\tfrac12x^\top Px+q^\top x\\
	\text{subject to}\quad &Ax+s=b,\quad s\in K,
	\end{aligned}
	\tag{P} \label{eq:conic-program}
\end{equation}
where $x\in\R^n$, $s\in\R^m$, $P\succeq0$, and $K$ is a closed convex proper cone.

Problems of this form arise in numerous applications of practical interest. Examples include constrained optimal control \cite{wang2010}, structural limit analysis \cite{makrodimopoulos2007} and stability analysis of fluid flows \cite{goulart2012}. The same formulation also covers support vector machines and lasso \cite{boyd2011}, image processing \cite{chambolle2011}, circuit design, and portfolio optimisation \cite{boyd2004}. Such problems also occur as subproblems in sequential convex programming methods for non-convex optimisation, including sequential quadratic programming (SQP) \cite[Ch.~18]{nocedal2006}, and as continuous relaxations in branch-and-bound methods for mixed-integer optimisation \cite{morrison2016}.

Interior-point methods \cite{wright1997} and first-order splitting methods \cite{boyd2011} (IPMs and FOMs respectively) are widely used to solve \eqref{eq:conic-program} numerically. Broadly speaking, FOMs have cheaper iterations but take more of them to reach modest accuracy, while IPMs can often reach high accuracy in few iterations \cite{stellato2020}. A principal computational trade-off is that FOMs can reuse a fixed KKT factorisation, whereas IPMs typically refactor a changing matrix to account for barrier curvature at the current iterate.

Nonetheless, the saddle matrices to be factored in, say, Clarabel \cite{goulart2026} and COSMO \cite{garstka2021} look remarkably similar. This resemblance motivates the central questions of this paper: under what metric choices and recovery rules does a FOM step agree with an IPM step, and what does this tell us about convergence and the cost of the linear solve in each?

Our main result is an explicit splitting--Newton correspondence, expressed through what we call the \textit{Moreau slack} $v$: an unconstrained vector from which one can recover multipliers and (ordinary) slacks satisfying exact or barrier-perturbed complementarity. This parameterisation eliminates complementarity from the equations to be solved, leaving a unified residual of primal and dual feasibility. Its roots recover KKT points when the barrier parameter $\mu$ is zero and central-path points when it is positive. Splitting applied to this residual uses cone projection at $\mu=0$ and a barrier proximal map at $\mu>0$; Newton's method at $\mu>0$ gives an ``implicit'' IPM.

At fixed $\mu>0$, matching the splitting projection metric to the barrier curvature at the current slack gives the same primal correction as implicit Newton with the same primal regularisation, while their Moreau-slack corrections differ by a factor of two. With this metric held fixed through the step, doubling only the Moreau-slack correction therefore gives the same finite iterate under a common step length and reconstruction.

We use this correspondence to examine three practical choices in algorithm design. Under FOM--IPM curvature matching, local analysis at fixed $\mu>0$ identifies when doubling improves (local) convergence: it removes a slack error mode but leaves the primal error mode unchanged. Metric choices that simplify the splitting solve act as regularisation in the matched Newton step; a two-variable example shows how making the linear solve cheap (\eg diagonal) can slow convergence. Finally, for orthant barrier subproblems, we derive a computable defect in the regularised Newton equations that measures the accuracy of a doubled splitting correction using a stored factorisation, providing a criterion for reusing the factors or refreshing the curvature and refactorising.

Our principal correspondence builds on established complementarity parameterisations and barrier methods. The Moreau slack uses the classical Minty parameterisation and Moreau decomposition \cite{minty1962,moreau1965}; a similar variable also appears in the ADMM analysis of Banjac et al.\ \cite{banjac2019}. Chua and Li \cite{chua2013} use barrier-based smoothing of cone projections to solve normal-map equations with Newton's method, while Arrizabalaga and Manchester \cite{arrizabalaga2026} recover complementary slacks and multipliers from an auxiliary variable in an implicit IPM for quadratic programming. Barrier-based splitting is exemplified by ABIP \cite{lin2021}, which uses ADMM to approximately solve successive barrier subproblems of a homogeneous self-dual formulation of a linear program.

The remainder of the paper is organised as follows. In \Cref{sec:tutorial} we derive the familiar linear systems underlying proximal ADMM and primal--dual IPMs; \Cref{sec:bridge} introduces the common residual and establishes the exact correspondence; and \Cref{sec:consequences} examines its implications for doubling corrections, simplifying solves, and reusing factorisations.

\subsection{Notation}
\label{sec:notation}

We use $\ip{u}{v}=u^\top v$ and $\|u\|=\|u\|_2$; for $Q\succ0$, $\|u\|_Q^2=u^\top Qu$. The symbols $\succeq$ and $\succ$ denote the positive semidefinite and positive definite matrix orders, respectively; vector inequalities are componentwise. The symbols $I$ and $\mathbf1$ denote the identity matrix and all-ones vector of appropriate size, and $\diag$ forms diagonal or block-diagonal matrices. We write $\rho(C)$ for the spectral radius of $C$ and $\lambda_{\max}(Q)$ for the largest eigenvalue of a symmetric matrix $Q$. A superscript $+$ denotes the next iterate.

For a cone $K\subseteq\R^m$, $\operatorname{int}K$ is its interior and $K^*=\{y\in\R^m:\ip{s}{y}\geq0\text{ for all }s\in K\}$ its dual; $\R_+^m$ is the nonnegative orthant. The indicator $\ind_C$ is zero on $C$ and $+\infty$ elsewhere. For a convex function $g$, $\partial g$ and $g^*$ denote its subdifferential and convex conjugate, respectively.

\section{Two familiar methods on the same problem}
\label{sec:tutorial}

In this section we apply to \eqref{eq:conic-program} two typical workhorse methods, one a FOM, the other an IPM, and demonstrate the suggestive similarity of the linear systems at their core.

Start by associating multiplier $y \in \R^m$ with \eqref{eq:conic-program} and define the primal and dual residuals
\begin{equation}
	r_p=Ax+s-b, \qquad r_d=Px+q+A^\top y.
 \label{eq:residuals}
\end{equation}

\subsection{Proximal ADMM}

Fix $M,R\succ0$. An unrelaxed proximal ADMM iteration is
\begin{equation}
 \begin{aligned}
 x^+&=\underset{\xi}{\arg\min}\bigl\{f(\xi)+\tfrac12\|\xi-x\|_M^2+\tfrac12\|A\xi+s-b+R^{-1}y\|_R^2\bigr\},\\
 s^+&=\Pi_K^R(b-Ax^+-R^{-1}y),\\
 y^+&=y+R(Ax^++s^+-b),
 \end{aligned}
 \label{eq:admm-native}
\end{equation}
where $\Pi_K^R$ is projection in the $R$-metric. The last two updates enforce cone complementarity after the first iteration. Writing $\Delta x=x^+-x$, the first-order condition of the primal subproblem gives
\begin{equation}
 (P+M+A^\top RA)\Delta x=-(r_d+A^\top Rr_p).
 \label{eq:split-schur}
\end{equation}
Introducing $\Delta\lambda=R(r_p+A\Delta x)$ rewrites this as the symmetric system
\begin{equation}
 \underbrace{\begin{bmatrix}P+M&A^\top\\A&-R^{-1}\end{bmatrix}}_{K_{\rm split}(M,R)}
 \begin{bmatrix}\Delta x\\\Delta\lambda\end{bmatrix}
 =-\begin{bmatrix}r_d\\r_p\end{bmatrix}.
 \label{eq:ksplit}
\end{equation}
Here $\Delta\lambda$ is an auxiliary solve variable: the finite multiplier change is $y^+-y=\Delta\lambda+R(s^+-s)$. For $M=\sigma I$ and $R = \diag(\rho_i)$ with positive parameters $\sigma$ and $\{\rho_i\}_i$, $K_{\rm split}$ is the quasi-definite matrix used in COSMO \cite{garstka2021} and OSQP \cite{stellato2020} (the latter being specialised to quadratic programs).

\subsection{Primal--dual interior-point methods}

Conventional (or ``explicit'') primal--dual IPMs retain (barrier-perturbed) complementarity as an equation alongside primal and dual feasibility. Let us consider the orthant $K=\R_+^m$ for a moment. At positive $s,y$, write $S=\diag(s)$, $Y=\diag(y)$, and $r_c=Sy-\mu\mathbf1$ for a target $\mu>0$. Linearising gives
\begin{equation}
 \begin{aligned}
 P\Delta x+A^\top\Delta y&=-r_d,\\
 A\Delta x+\Delta s&=-r_p,\\
 Y\Delta s+S\Delta y&=-r_c.
 \end{aligned}
\end{equation}
Eliminating $\Delta s$ yields~\cite{vandenberghe2010}
\begin{equation}
 \begin{bmatrix}P&A^\top\\A&-Y^{-1}S\end{bmatrix}
 \begin{bmatrix}\Delta x\\\Delta y\end{bmatrix}
 =\begin{bmatrix}-r_d\\-r_p+Y^{-1}r_c\end{bmatrix}.
 \label{eq:pdipm}
\end{equation}
The inverse constraint weight is now $Y^{-1}S$, determined by the current slack and multiplier. Along an orthant central path, the weight itself is $y_i/s_i=\mu/s_i^2$: approaching a strictly complementary solution, it grows on active constraints and shrinks on inactive ones. In the ADMM system~\eqref{eq:ksplit}, the analogous weights are the penalties $\rho_i$: choosing larger penalties for constraints judged active and smaller ones for constraints judged inactive therefore mimics this IPM weighting, and the idea is used to varying degrees in such solvers as OSQP \cite{stellato2020}, COSMO \cite{garstka2021}, and SuperADMM \cite{verheijen2025}.

For symmetric cones in general, the coefficient block $Y^{-1}S$ is usually replaced by $W^\top W$, where $W$ is the Nesterov--Todd scaling matrix satisfying $Wy=W^{-\top}s$~\cite[Secs.~4--5]{vandenberghe2010}. In the orthant case, $W=\diag(\sqrt{s_i/y_i})$. Thus $W^\top W$ occupies the same block as $R^{-1}$ in~\eqref{eq:ksplit}. This saddle structure underlies Clarabel's linear algebra~\cite{goulart2026}.

\section{A common residual and an exact correspondence}
\label{sec:bridge}

\subsection{Eliminating complementarity}

To compare splitting and Newton on the same equations, we treat projection and barrier complementarity in one family parameterised by $\mu \geq 0$.

Let $F$ be a logarithmically homogeneous self-concordant barrier for $K$, extended by $+\infty$ outside $\operatorname{int}K$, and define $\psi_\mu=\mu F$ for $\mu>0$ and $\psi_0=\ind_K$ (the indicator function of set $K$). The choice at zero is natural in the sense that $\mu F$ epi-converges to $\ind_K$ as $\mu\downarrow0$~\cite[Thm.~7.17]{rockafellar1998}.

The primal and dual residuals in~\eqref{eq:residuals}, together with a parameterised complementarity relation, give the central path conditions
\begin{equation}
 r_d=0,\qquad r_p=0,\qquad -y\in\partial\psi_\mu(s).
 \label{eq:kkt}
\end{equation}
At $\mu=0$, the last relation is $s\in K$, $y\in K^*$, and $\ip{s}{y}=0$, encoding KKT points of \eqref{eq:conic-program}. For $\mu>0$, it is $y=-\mu\nabla F(s)$, encoding central path points. In the orthant example, with the standard barrier function $F(s)=-\sum_i\log s_i$, this is the familiar perturbed complementarity condition $s_i y_i=\mu$ with $s,y>0$.

We now seek to represent the slack and multiplier by a single unconstrained vector, so that complementarity is satisfied automatically. Fix a positive-definite metric $R\succ0$, and define the weighted proximal map
\begin{equation}
 p_{\mu,R}(v)=\underset{s}{\arg\min}\left\{
 \psi_\mu(s)+\tfrac12(s-v)^\top R(s-v)\right\}.
 \label{eq:prox}
\end{equation}
At $\mu=0$, $p_{0,R}$ is projection onto $K$ in the $R$-metric, $\Pi_K^R$. For each fixed $R\succ0$ and $v$, epi-convergence also gives $p_{\mu,R}(v)\to\Pi_K^R(v)$ as $\mu\downarrow0$~\cite[Thm.~12.35]{rockafellar1998}.

We call $v$ the \emph{Moreau slack}. It encodes a slack--multiplier pair satisfying the prescribed complementarity relation in one unconstrained vector, from which both variables are recovered by
\begin{equation}
 \boxed{s=p_{\mu,R}(v),\qquad y=R(s-v).}
 \label{eq:minty}
\end{equation}
The proximal optimality condition $0\in\partial\psi_\mu(s)+R(s-v)$ shows that~\eqref{eq:minty} is a one-to-one parameterisation of all slack--multiplier pairs satisfying the complementarity relation in~\eqref{eq:kkt}, with inverse $v=s-R^{-1}y$. This is the classical Minty parameterisation~\cite{minty1962}, here written in the metric $R$.

For the nonnegative orthant $K=\R_+^m$, with $F(s)=-\sum_i\log s_i$ and $R=\diag(r_i)$, $r_i>0$, recovery through~\eqref{eq:minty} takes the explicit form
\begin{equation}
 s_i=\frac{v_i+\sqrt{v_i^2+4\mu/r_i}}{2},
 \qquad y_i=r_i(s_i-v_i).
 \label{eq:orthant-recovery}
\end{equation}
At $\mu=0$, this reduces to $s_i=\max\{v_i,0\}$; for $\mu>0$, both $s_i$ and $y_i$ are strictly positive and satisfy $s_i y_i=\mu$. Thus the same unconstrained coordinate $v$ supports both projection and barrier recovery. Related barrier-based constructions appear in Chua and Li~\cite{chua2013} and Arrizabalaga and Manchester~\cite{arrizabalaga2026}. For $\mu=0$ and $R=\rho I$, Banjac et al.'s variable $v^k=z^k+\rho^{-1}y^k$~\cite[Sec.~4.1, eq.~(13)]{banjac2019} plays the same role, corresponding to $b-v$ when $z=b-s$.

Our name ``Moreau slack'' reflects Moreau's decomposition~\cite{moreau1965}: once $s=p_{\mu,R}(v)$ is recovered, the remainder gives the multiplier through the convex conjugate $\psi_\mu^*$,
\begin{equation}
	\prox_{\psi_\mu^*}^{R^{-1}}(Rv)
	=R\bigl(v-p_{\mu,R}(v)\bigr)
	=-y \quad\implies\quad v=\underbrace{\prox_{\psi_\mu}^{R}(v)}_{s}+R^{-1}\underbrace{\prox_{\psi_\mu^*}^{R^{-1}}(Rv)}_{-y}.
\end{equation}
Here the superscripts denote the proximal metrics. The identity holds for every $\mu\geq0$; at $\mu=0$, the proximal maps are projections to the cones $K$ and $-K^*$ (the form used in the ADMM analysis of \cite{banjac2019}).

Substituting~\eqref{eq:minty} into the remaining feasibility equations leaves the residual
\begin{equation}
 \boxed{\res_{\mu,R}(x,v)=
 \begin{bmatrix}
 Px+q+A^\top R\bigl(p_{\mu,R}(v)-v\bigr)\\
 Ax+p_{\mu,R}(v)-b
 \end{bmatrix}.}
 \label{eq:residual}
\end{equation}
Its roots reconstruct central-path points when $\mu>0$ and KKT points when $\mu=0$.

We now show how a splitting iteration and Newton's method can both be used to annihilate this residual. At $\mu=0$, splitting gives a projection-based FOM, such as ADMM in its familiar forms. At $\mu>0$, it instead acts on the barrier subproblem
\begin{equation}
 \min_{x,s}\ f(x)+\mu F(s)\quad\text{subject to}\quad Ax+s=b.
\end{equation}
We call this \emph{barrier splitting}: a first-order method for a fixed barrier subproblem, with cone projection replaced by a barrier proximal map. ABIP~\cite{lin2021}, for example, uses ADMM to approximately solve successive barrier subproblems of a homogeneous self-dual formulation of a linear program.

Here we focus on the iteration at \emph{fixed} $\mu$; beyond this, a complete path-following method (along the lines of ABIP) would also specify how to decrease $\mu$ and how accurately to solve each subproblem.

\subsection{Two ways of solving the feasibility equations}

\subsubsection{Splitting}

Let us start with a splitting iteration to solve $\res_{\mu,R}(x,v) = 0$. Write the objective with its affine constraint as
\begin{equation}
 \phi(x,s)=f(x)+\ind_{\{Ax+s=b\}}(x,s).
\end{equation}
Apply Douglas--Rachford splitting (DRS) to $\phi(x,s)+\psi_\mu(s)$ in the product metric $\diag(M,R)$, with $M\succ0$. Using the ordering $\prox_{\psi_\mu}$ then $\prox_\phi$, one iteration in $u=(x,v)$ is
\begin{equation}
\begin{aligned}
 s&=p_{\mu,R}(v),\\
 (\widetilde x,\widetilde s)&=\prox_\phi^{\diag(M,R)}(x,2s-v),\\
 (x^+,v^+)&=(x,v)+\alpha(\widetilde x-x,\widetilde s-s),
 \quad 0<\alpha<2.
\end{aligned}
\end{equation}
The second proximal step is the solution of an equality-constrained QP. Writing $\Delta x=\widetilde x-x$ and $\Delta v=\widetilde s-s$, its optimality conditions give
\begin{equation}
 \underbrace{\begin{bmatrix}
 P+M&-A^\top R\\ A&I
 \end{bmatrix}}_{L_{M,R}}
 \begin{bmatrix}\Delta x\\\Delta v\end{bmatrix}
 =-\res_{\mu,R}(x,v).
 \label{eq:dr-correction}
\end{equation}
Hence DRS is precisely
\begin{equation}
	u^+=u-\alpha L_{M,R}^{-1}\res_{\mu,R}(u).
\end{equation}
Setting $\Delta\lambda=-R\Delta v$ recovers the symmetric splitting system~\eqref{eq:ksplit}, and eliminating $\Delta v$ gives~\eqref{eq:split-schur}. For $\alpha=1$, the new coordinate is $v^+=b-Ax^+-R^{-1}y$. Thus, at $\mu=0$ reconstruction through \eqref{eq:minty} gives the ADMM updates~\eqref{eq:admm-native} once the current slack and multiplier satisfy complementarity. (Positive $\mu$ smoothes the projection into the barrier prox.)

The close relationships between ADMM and DRS are well known in the convex optimisation literature~\cite{banjac2019}. Similar connections extend to other FOMs: the primal--dual hybrid gradient method (PDHG) admits a DRS representation~\cite{oconnor2020}; this as well as preconditioned ADMM~\cite{bredies2017} and several other splitting schemes~\cite{bredies2022} can be viewed as preconditioned proximal point methods.

At each fixed $\mu\geq0$, exact DRS with fixed positive-definite metrics and constant $0<\alpha<2$ converges to a fixed point whose reconstruction satisfies~\eqref{eq:kkt}, provided such a point exists~\cite[Thm.~7]{eckstein1992}.

\subsubsection{Newton's method}

Now consider solving $\res_{\mu,R}(x,v)=0$ using Newton's method, with $\mu>0$ and $R$ fixed.
Write $s=p_{\mu,R}(v)$ and $B=\mu\nabla^2F(s)$. Differentiating \eqref{eq:prox}'s proximal optimality condition $\mu\nabla F(s)+R(s-v)=0$ with respect to $v$ gives
\begin{equation}
	D:=Dp_{\mu,R}(v)=(R+B)^{-1}R.
\end{equation}
The Jacobian of the common residual with respect to $u=(x,v)$, evaluated at the current iterate, is then
\begin{equation}
 J_{\mu,R}:=D_u\res_{\mu,R}=\begin{bmatrix}
 P&A^\top R(D-I)\\ A&D
 \end{bmatrix}.
 \label{eq:jac}
\end{equation}
The corresponding linearised changes are $\Delta s=D\Delta v$ and $\Delta y=-B\Delta s$.\footnote{For the latter, note that \eqref{eq:minty} implies $\Delta y = R (D - I) \Delta v$, and $R(D-I)=-BD$.} The system is singular precisely when the null spaces of $P$ and $A$ have a non-trivial intersection. In practice one often uses regularisation to ensure definiteness~\cite{goulart2026}. It follows from the above that an implicit Newton correction with primal regulariser $M$ solves
\begin{equation}
 \bigl(J_{\mu,R}+\diag(M,0)\bigr)\Delta u_{\rm N}=-\res_{\mu,R}(u).
 \label{eq:newton-correction}
\end{equation}

\subsection{Matching curvature and recovering the iterate}

When the metric $R$ in barrier splitting matches the current barrier Hessian $B$, we have $D=I/2$. Thus, recovery through~\eqref{eq:minty} changes the slack $s$ only half as far as $v$ to first order, accounting for the factor of two below.\footnote{No finite choice of $R$ gives $D=I$ when $B\ne0$, since $(R+B)D=R$ would then imply $B=0$. One can approach $D=I$ by taking $R=tB$ and letting $t\to\infty$, for which $D=\frac{t}{t+1}I$.}

\begin{proposition}[Exact curvature matching]
\label{prop:exact-matching}
Fix $\mu>0$ and a physical state $(x,s,y)$ with $s\in\operatorname{int}K$ and $y=-\mu\nabla F(s)$. Choose $M\succ0$, set $R=B=\mu\nabla^2F(s)$, and encode the state by $v=s-R^{-1}y$. Use the same $M$ in barrier splitting~\eqref{eq:dr-correction} and primal-regularised Newton~\eqref{eq:newton-correction}, holding $\mu,R$ fixed through the solve and reconstruction.

The corresponding corrections satisfy
\begin{equation}
 \boxed{\Delta x_{\rm N}=\Delta x_{\rm DR},\qquad
 \Delta v_{\rm N}=2\Delta v_{\rm DR}.}
 \label{eq:directions}
\end{equation}
Doubling only $\Delta v_{\rm DR}$ and using the same step length and nonlinear recovery~\eqref{eq:minty} gives the same finite iterate as implicit Newton.
\end{proposition}

\emph{Proof.} Substituting $D=I/2$ into~\eqref{eq:jac} gives
\begin{equation}
 J_{\mu,R}+\begin{bmatrix}M&0\\0&0\end{bmatrix}
 =L_{M,R}\begin{bmatrix}I&0\\0&\tfrac12I\end{bmatrix}.
 \label{eq:match}
\end{equation}
Both systems are nonsingular because their common primal Schur complement $P+M+A^\top BA$ is positive definite. Comparing their solutions proves~\eqref{eq:directions}; the same recovery gives identical iterates. \hfill$\square$

To repeat the correspondence at the new physical state, rematch $R=B(s)$ and re-encode $v=s-R^{-1}y$.

On the orthant complementarity graph, $B^{-1}=\diag(s_i/y_i)$, so the implicit Newton tangents agree with those of the conventional primal--dual system~\eqref{eq:pdipm} under the same primal regularisation $M$. Their finite updates generally differ: the additive updates $s^+=s+\alpha\Delta s$ and $y^+=y+\alpha\Delta y$ give
\begin{equation}
 (s_i+\alpha\Delta s_i)(y_i+\alpha\Delta y_i)
 =\mu-\alpha^2B_{ii}(\Delta s_i)^2,
\end{equation}
whereas the matched step recovers $(s^+,y^+)$ from $v^+=v+2\alpha\Delta v_{\rm DR}$ via~\eqref{eq:minty}, preserving perturbed complementarity exactly.

\section{Consequences for algorithm design}
\label{sec:consequences}

Our correspondence above suggests three immediate questions: whether to double the Moreau-slack correction in barrier splitting; how far to simplify the linear solve in splitting or Newton methods through tailored choices of $M$; and when to reuse a factorisation with curvature from a previous iterate. We first identify the gain from doubling; then examine how the solve limits that gain; and finally use the correspondence to guide factorisation reuse.

\subsection{What does doubling $\Delta v$ buy?}
\label{sec:doubling}

Doubling gives the barrier splitting step a Newton interpretation in \eqref{eq:directions}, but does it make the iteration faster? We answer this locally by freezing the metrics at a barrier solution and linearising the iterations~\cite{liang2017}. With frozen metrics, doubled barrier splitting and implicit Newton generally agree only to first order.

\begin{proposition}[Local effect of doubling]
\label{prop:common-local}
Fix $\mu>0$ and a barrier solution $(x_*,s_*,y_*)$, with $F$ of class $C^3$ near $s_*$. Let $B_*=\mu\nabla^2F(s_*)$ and suppose the unregularised barrier Schur complement is definite, \ie $H_0=P+A^\top B_*A\succ0$. For any $M_*\succ0$, set $H_*=H_0+M_*$ and define
\begin{equation}
 C=H_*^{-1}M_*,\qquad c=\rho(C)<1.
 \label{eq:local-primal-factor}
\end{equation}

Freeze $R=B_*$ and $M=M_*$. Unrelaxed barrier DRS and the (``doubled'') iteration obtained by doubling only its correction in $v$, both with recovery~\eqref{eq:minty}, are locally linearly convergent. To first order, both methods give the same primal error,
\begin{equation}
 x^+-x_* = C(x-x_*).
\end{equation}
Their slack errors satisfy, to first order,
\begin{equation}
 s^+-s_* = \begin{cases}
 \tfrac12(s-s_*)+\tfrac12\bigl[-A(x^+-x_*)\bigr],&\text{barrier DRS},\\
 -A(x^+-x_*),&\text{doubled}.
 \end{cases}
 \label{eq:local-slack-errors}
\end{equation}
Their local spectral radii are respectively
\begin{equation}
 \rho_{\rm DRS}=\max\{c,\tfrac12\},\qquad
 \rho_{\rm doubled}=c.
 \label{eq:local-spectral-radii}
\end{equation}
Unit-step, $M_*$-regularised implicit Newton with the same recovery has the same first-order recurrence as doubled barrier splitting.
\end{proposition}

\emph{Proof.} Write $e_x=x-x_*$, $e_s=s-s_*$ and $e=(e_x,e_s)$. Recovery through~\eqref{eq:minty} enforces $y=-\mu\nabla F(s)$, so Taylor expansion at $s_*$ gives $y-y_*=-B_*e_s+o(\|e\|)$. Using~\eqref{eq:residuals} and the vanishing residuals at the solution, the slack terms cancel to first order:
\begin{equation}
 \begin{aligned}
 r_d+A^\top B_*r_p
 &=Pe_x-A^\top B_*e_s+A^\top B_*(Ae_x+e_s)+o(\|e\|)\\
 &=H_0e_x+o(\|e\|).
 \end{aligned}
\end{equation}
Thus~\eqref{eq:split-schur} gives $\Delta x=-H_*^{-1}H_0e_x+o(\|e\|)$, and hence $e_x^+=Ce_x+o(\|e\|)$ for both splitting methods.

The second row of~\eqref{eq:dr-correction} gives $\Delta v_{\rm DR}=-r_p-A\Delta x=-e_s-ACe_x+o(\|e\|)$. At the solution, $Dp_{\mu,R}(v_*)=(R+B_*)^{-1}R=I/2$. Consequently, reconstructing from $v^+=v+\Delta v_{\rm DR}$ or $v^+=v+2\Delta v_{\rm DR}$ gives $s^+-s=\tfrac12\Delta v_{\rm DR}+o(\|e\|)$ or $s^+-s=\Delta v_{\rm DR}+o(\|e\|)$, respectively. Adding these changes to $e_s$ yields~\eqref{eq:local-slack-errors}.

The first-order error maps on $(e_x,e_s)$ for DRS and doubled splitting are therefore
\begin{equation}
 T_{\rm DRS}=\begin{bmatrix}C&0\\-\tfrac12AC&\tfrac12I\end{bmatrix},
 \qquad
 T_{\rm doubled}=\begin{bmatrix}C&0\\-AC&0\end{bmatrix},
\end{equation}
and their block-triangular form gives the stated spectra. Regularised Newton has the same first-order error map as doubled splitting: its coefficient matrix at the solution satisfies~\eqref{eq:match}, and its variation contributes only higher-order terms because the residual vanishes there.

Since $0\prec M_*\prec H_*$, we have $c<1$, so both first-order error maps have spectral radius strictly below one. Local linear convergence of the corresponding differentiable iterations then follows from Ostrowski's fixed-point theorem~\cite[Thm.~10.1.3]{ortega2000}. \hfill$\square$

\paragraph{Interpreting the local rates.}
Since $Ax_*+s_*=b$, the slack that restores feasibility at $x^+$ is $b-Ax^+$, whose deviation from $s_*$ is $-A(x^+-x_*)$. To first order, barrier DRS moves the slack halfway from $s$ to $b-Ax^+$, while the doubled correction takes it all the way.

Doubling removes the independent slack error mode: to first order, a perturbation with $e_x=0$ is halved by barrier DRS and eliminated by the doubled correction. The remaining slack error, $e_s^+=-ACe_x$, is induced by primal error, so both are governed by $c$.

The common primal factor $c$ measures the splitting metric $M_*$ relative to the barrier subproblem's primal curvature $H_0$: in a nonzero direction $z$, their ratio is $z^\top M_*z/(z^\top H_0z)$. Since $M_*z=\lambda H_0z$ implies $Cz=\lambda z/(1+\lambda)$, the largest such ratio determines the slowest primal mode:
\begin{equation}
 \theta=\lambda_{\max}\!\left(H_0^{-1/2}M_*H_0^{-1/2}\right),
 \qquad c=\frac{\theta}{1+\theta}.
\end{equation}
Doubling therefore improves the local factor over unrelaxed DRS exactly when $c<1/2$, or equivalently when $M_*\prec H_0$; otherwise the primal solve remains the bottleneck. This requires regularisation to be smaller than primal curvature in every direction; for $M_*=\sigma I$, it means $\sigma<\lambda_{\min}(H_0)$. Also, note that ordinary scalar over-relaxation can give a smaller factor still: for $M_*=(2/3)H_0$, ordinary DRS with $\alpha=20/11$ has factor $1/11$, compared with $2/5$ after doubling only the Moreau-slack correction.

The assumption $H_0\succ0$ permits $M=0$ in both~\eqref{eq:dr-correction} and~\eqref{eq:newton-correction}: the systems remain nonsingular near the solution, and unregularised Newton and the doubled correction with $R=B_*$ both have $C=0$ and an $O(\|e\|^2)$ error recurrence under the stated $C^3$ assumption, giving local quadratic convergence.

\paragraph{Limits of doubling.}
Positive $\mu$ alone does not guarantee that doubling helps. For $\min_x\tfrac12x^2$ subject to $-1\le x\le1$, take $A=(1,-1)^\top$, $b=\mathbf1$, $M=1$, $R=I$, and the standard logarithmic barrier. States with $x=0$ and $v=t\mathbf1$ remain of this form, with $t^+=t+\beta[1-p_\mu(t)]$, where $p_\mu(t)=(t+\sqrt{t^2+4\mu})/2$ is the scalar $R=I$ specialisation of~\eqref{eq:orthant-recovery} and $\beta=1$ or $2$. At $t_*=1-\mu$, the local factors along this trajectory are
\begin{equation}
 \left|\frac{dt^+}{dt}(t_*)\right|=\begin{cases}
 \dfrac{\mu}{1+\mu},&\text{ordinary DRS},\\[4pt]
 \dfrac{|1-\mu|}{1+\mu},&\text{doubled}.
 \end{cases}
\end{equation}
For $\mu=1/10$, doubling worsens the factor from $1/11$ to $9/11$. At $\mu=0$ and $0<t<2$, the updates become $t^+=1$ and $t^+=2-t$, giving an exact doubled two-cycle (\ie destroying convergence) unless $t=1$.

For $\mu>0$, we next examine how simplifying the primal solve through our choice of $M$ changes the remaining factor $c$.

\subsection{What is the price of simplifying the primal solve?}
\label{sec:solve-structure}

We can use a tailored choice of $M$ to remove coupling directions and cheapen the required primal solve, as we demonstrate below. As usual, there is no free lunch: this can slow the primal mode considerably.

At matching, write $H_{\rm N}=P+A^\top BA$. We can make the resulting primal solve simple (say, diagonal) by choosing a simple positive-definite majorant $H$, so that $H\succ H_{\rm N}$ and the primal regulariser is
\begin{equation}
	M=H-H_{\rm N}\succ0.
 \label{eq:structured-majorant}
\end{equation}
The idea is not new: non-diagonal regularisation in an IPM was similarly used in~\cite{pougkakiotis2019}; in splitting, $M=T^{-1}-A^\top RA\succ0$ with diagonal positive definite $T$ is the PDHG metric choice, leaving primal coefficient $P+T^{-1}$~\cite{chambolle2011,pock2011,oconnor2020} (diagonal if $P$ is diagonal). Going beyond PDHG, one can clearly choose a majorant that eliminates not only the constraint Gramian matrix but also the off-diagonal entries of the Hessian, resulting in a diagonal solve even for general Hessian matrices.

Suppose the chosen majorant $H=H(x,s)$ depends continuously on the current iterate. Rematching $R=B(s)$ and re-encoding $v$ at each step, with $M=H-H_{\rm N}$ and doubled recovery, gives the local convergence factor $c$ from Proposition~\ref{prop:common-local}, evaluated with $H_*=H(x_*,s_*)$ and $M_*=H_*-H_0$. As in the proof, variation of the solve matrix contributes only higher-order terms because the residual vanishes at the solution.

We use the following example to illustrate the convergence costs of removing coupling:
\begin{equation}
 \min_{x\in\R^2}\ \tfrac12\|x-\mathbf1\|^2
 \quad\text{subject to}\quad x_1+x_2\le1.
 \label{eq:paired-example}
\end{equation}
For ordinary PDHG, choose $R=\rho\ge1$ and $T=(1+2\rho)^{-1}I$, so $M\succ0$ and the primal coefficient is $(2+2\rho)I$. At a positive-barrier solution with $B=\rho$, choose the same diagonal coefficient as the majorant:
\begin{equation}
 H_{\rm N}=I+\rho\begin{bmatrix}1&1\\1&1\end{bmatrix},
 \qquad H_{\rm diag}=(2+2\rho)I.
\end{equation}
This is a valid majorant: $M=H_{\rm diag}-H_{\rm N}=(1+2\rho)I-\rho\mathbf1\mathbf1^\top$ has eigenvalues $1$ and $1+2\rho$.

For each method, take the same primal error $e_x=x-x_*=t(1,-1)^\top$ about its respective solution, leaving its slack and multiplier unchanged. Then $r_p=0$, $r_d=e_x$, and $\Delta v_{\rm DR}=0$. Retaining the coupling in the interior solve with $H_{\rm pair}=H_{\rm N}+I$ gives the comparison
\begin{equation}
 \begin{aligned}
 e_{x,\rm PDHG}^+=e_{x,\rm diag}^+
   &=\left(1-\frac1{2+2\rho}\right)e_x,\\
 e_{x,\rm pair}^+&=\tfrac12 e_x.
 \end{aligned}
 \label{eq:paired-tangent-factors}
\end{equation}
Constraint curvature acts normally to the constraint; diagonalisation spreads it into the tangent direction, making the first factor tend to one as $\rho$ grows. For the matched barrier step, $\rho=B=y^2/\mu\sim1/(4\mu)$ as $\mu\downarrow0$, since $sy=\mu$ and $y\to1/2$ along the barrier path (ordinary PDHG may keep $\rho$ fixed, of course). Retaining the coupling prevents the slowdown.


\subsection{When should curvature be refreshed?}
\label{sec:factorisation-reuse}

Exact matching recovers Newton, but we can assess whether a stored metric remains adequate by measuring the current correction's regularised Newton-equation defect. Quasi-Newton IPMs also exploit factorisation reuse~\cite{gondzio2019}. Note that the orthant specialisation below makes the curvature comparison componentwise (hence cheap), though the defect identity itself holds for general cones/barriers.

Fix $\mu>0$, $M\succ0$, and an orthant state with $s>0$, $y_i=\mu/s_i$; we also allow $M=0$ when $P+A^\top BA\succ0$. For a stored diagonal $R\succ0$, encode $v=s-R^{-1}y$ and compute $(\Delta x,\Delta v_{\rm DR})$ from~\eqref{eq:dr-correction}. The doubled correction $\widehat{\Delta u}=(\Delta x,2\Delta v_{\rm DR})$ has regularised Newton-equation defect
\begin{equation}
 \begin{aligned}
 e&:=\bigl(J_{\mu,R}+\diag(M,0)\bigr)\widehat{\Delta u}+\res_{\mu,R}\\
  &=\begin{bmatrix}A^\top Rz\\z\end{bmatrix},\qquad
 z_i=\frac{R_{ii}-B_{ii}}{R_{ii}+B_{ii}}(\Delta v_{\rm DR})_i.
 \end{aligned}
 \label{eq:refresh-defect}
\end{equation}
Indeed, subtracting~\eqref{eq:dr-correction} from the linearised equations~\eqref{eq:jac} leaves the blocks $A^\top R(2D-I)\Delta v_{\rm DR}$ and $(2D-I)\Delta v_{\rm DR}$. In general, $z=(R+B)^{-1}(R-B)\Delta v_{\rm DR}$, which reduces to~\eqref{eq:refresh-defect} for diagonal $R$ and $B$.

Using inexact-Newton residual control~\cite{dembo1982}, an illustrative curvature refresh rule would retain the factors when
\begin{equation}
 \|e\|_\infty\leq\eta_{\rm N}\|\res_{\mu,R}\|_\infty,
 \qquad 0<\eta_{\rm N}<1.
\end{equation}
Otherwise, set $R=B$, re-encode $v$, and recompute the correction at the same physical state, keeping $M$ fixed. The refreshed defect is zero.
The unregularised Newton defect is $e-(M\Delta x,0)^\top$, and convergence under changing metrics would require further analysis still.

The following example shows why curvature mismatch alone does not determine whether to refresh, with or without primal regularisation. Return to~\eqref{eq:paired-example} at $\mu=1$, and take $M=mI$, $m\in\{0,1\}$, with stored $R=3$. At $x=t(1,-1)^\top$, $t\ne0$, and $s=y=1$, we have $B=1$, $\Delta v_{\rm DR}=0$, and $e=0$; refreshing gives the identical step. At $x=0$, $s=1/3$, $y=3$, we instead have $B=9$, $\Delta x=0$, $\Delta v_{\rm DR}=2/3$, and
\begin{equation}
 e=-(1,1,1/3)^\top,\qquad \res=(2,2,-2/3)^\top.
\end{equation}
Both states have $|(R-B)/(R+B)|=1/2$, but their relative defects are $0$ and $1/2$: $\eta_{\rm N}=1/8$ retains the first correction and refreshes the second for either choice of $m$. Note that for $m=0$, the test controls the unregularised Newton equations; for $m=1$, it controls their regularised counterpart. So the correspondence shows us how to check whether stored splitting factors still provide the desired Newton-equation accuracy before paying for a new factorisation.

\section{Conclusion}

The Moreau slack makes the resemblance between splitting and interior-point methods precise through a common residual. At fixed positive barrier parameter, curvature matching gives the same primal correction in barrier splitting and implicit Newton with the same primal regularisation; doubling the Moreau-slack correction then gives the same finite iterate under a common step length and nonlinear recovery.

This correspondence helps explain some algorithms' observed behaviour. With metrics frozen at a barrier solution, doubling the Moreau slack correction removes a local slack error mode, while regularisation determines the primal error that remains. Simplifying the solve can make this error decay arbitrarily slowly by discarding useful coupling. For orthant problems, the defect in the regularised Newton equations measures how reusing a factorisation affects the current correction. Together, these observations suggest adaptive barrier splitting methods that retain useful coupling and reuse factors while the corrections remain verifiably Newton-like.

\begingroup
\small
\bibliographystyle{IEEEtran}
\bibliography{paper/references}
\endgroup

\end{document}